\documentclass[11pt]{article}
\usepackage{amsmath,amsthm,amssymb,epsfig}

\newcommand\hyp{\mbox{\bf H}}
\newcommand\A{\mbox{\bf H}^2}
\newcommand\R{\mbox{\bf R}}
\newcommand\Z{\mbox{\bf Z}}
\newcommand\N{\mbox{\bf N}}
\newcommand\Cp{\mbox{\bf C}}
\newcommand\Q{\mbox{\bf Q}}
\newcommand\Qp{\mbox{\bf Q}_p}
\newcommand\Zp{\mbox{\bf Z}_p}
\newcommand\HT{\A \times T_p}
\newcommand\HTq{\A \times T_q}

\newcommand\F{\Z [\frac{1}{p}]}
\newcommand\Fq{\Z [\frac{1}{q}]}
\newcommand\PZp{PSL_2( \F)}
\newcommand\PZq{PSL_2( \Fq)}

\newcommand\PR{PSL_2( \R)}
\newcommand\G{PSL_2(\R) \times PSL_2( \Qp)}
\newcommand\PQp{PSL_2(\Qp)}

\newcommand\si{\sigma_{\infty}}
\newcommand\so{\sigma_0}

\newcommand\sa{\sigma_{\alpha}}
\newcommand\dq{\Delta_{Q}}
\newcommand\dz{\Delta}
\newcommand\cz{ \frac{1}{M}\Delta}
\newcommand\kz{ \frac{1}{k}\Delta}
\newcommand\qz{ \frac{1}{q}\Delta}

\newtheorem{theorem}{Theorem}[section]
\newtheorem{proposition}[theorem]{Proposition}
\newtheorem{lemma}[theorem]{Lemma}
\newtheorem{corollary}[theorem]{Corollary}

\input{epsf.tex}

\title{Quasi-isometric rigidity for $PSL_2({\bf Z}[ \frac{1}{p} ]$}
\author{Jennifer Taback}

\begin{document}

\maketitle

\abstract

We prove that $\PZp$ gives the first example of groups
which are not quasi-isometric to each other but have the same
quasi-isometry group.  Namely, $\PZp$ and $\PZq$ are not
quasi-isometric unless $p=q$, and, independent of $p$, the
quasi-isometry group of $\PZp$ is $PSL_2(\Q)$.  
In addition, we characterize $\PZp$ uniquely among all finitely generated groups
by its quasi-isometry type.

\section{Introduction}
\label{sec:results}
\pagenumbering{arabic}

Combining the work of many people yields a complete quasi-isometry
classification of irreducible lattices in semisimple Lie groups (see
\cite{F} for an overview of these results).
One of the first general results in this classification is the
complete description, up to quasi-isometry, of all nonuniform lattices
$\Lambda$ 
in semisimple Lie groups of rank $1$, proved by R. Schwartz \cite{S2}.  
He shows that every quasi-isometry of such a lattice $\Lambda$ is equivalent to a
unique commensurator of $\Lambda$.  (A {\it commensurator} of $\Lambda \subset G$ is an element $g \in G$ so that $g \Lambda g^{-1} \cap
\Lambda$ has finite index in $\Lambda$.)  
We will call this result {\em commensurator rigidity}, although it is a
different notion than the commensurator rigidity of Margulis.    
In \cite{FS} it was conjectured that commensurator rigidity, or at least a slightly
weaker statement, ``quasi-isometric iff commensurable,'' should apply
to nonuniform lattices in a wide class of Lie groups.
Here we prove that both of these statements are true for $\PZp$.

In a different direction, B. Farb and L. Mosher proved analogous
quasi-isometric rigidity results for the the solvable Baumslag-Solitar
groups.  
These groups are given by the
presentation 
$$ BS(1,n) = < a,b | aba^{-1} = b^{n} >$$
and are not lattices in any Lie group.

The group $\PZp$ is a nonuniform (i.e. non cocompact) lattice in $\G$, analogous to the
classical Hilbert modular group $PSL_2({\mathcal O}_d)$ in $PSL_2(\R) \times PSL_2(\R)$.  
It is also a basic example of an $S$-arithmetic group.  
The proofs of Theorems A, B and C (stated below) combine techniques
from the 
two types of quasi-isometric rigidity results mentioned above.  
When we construct a space $\Omega_p$ on which $\PZp$ acts properly
discontinuously and cocompactly by isometries, we see that the horospheres forming
the boundary components of $\Omega_p$ carry the geometry of the group
$BS(1,p)$.  
In this way the results of \cite{FM} play a role in the quasi-isometric rigidity of
$\PZp$.

\subsection{Statement of Results}

In this paper we prove the following  quasi-isometric rigidity results for the
finitely generated groups $\PZp$, where $p$ is a prime.  
Theorem A may be viewed as a strengthening of strong (Mostow) rigidity
for $\PZp$.  \cite{M}

\medskip
\noindent
{\bf Theorem A (Main Theorem)}.
{\it 
Every quasi-isometry of $\PZp$ 
is equivalent to a commensurator of
$\PZp$.  Hence the natural map
$$ Comm(PSL_2(\Z[ {\hbox{$\frac{1}{p}$}}])) \rightarrow QI(PSL_2(\Z[ {\hbox{$\frac{1}{p}$}}]))$$
is an isomorphism.}
\medskip

Since for any prime $p$ the commensurator group of $\PZp$ in $\G$ is $PSL_2(\Q)$, the
quasi-isometry group is also $PSL_2(\Q)$. (See \S \ref{sec:algebraic} for a model of $PSL_2$
as an algebraic group.)  
Thus we cannot distinguish the quasi-isometry classes of these groups
via their quasi-isometry groups.  
However, using a result of B. Farb and L. Mosher \cite{FM} (Theorem
\ref{thm:BS:qi} below), we are able to prove the
following.  

\medskip
\noindent
{\bf Theorem B (Quasi-isometric iff commensurable).}
{\it 
Let $p$ and $q$ be primes.  
Then $\PZp$ and $\PZq$ are quasi-isometric if and only if they are
commensurable, which occurs only when $p = q$.}
\medskip

Theorems A and B together give the first example of groups which have
the same quasi-isometry group but are not quasi-isometric.

The following theorem characterizes $\PZp$ uniquely among all
finitely generated groups by its quasi-isometry type.

\medskip
\noindent
{\bf Theorem C (Quasi-isometry characterization).}
{\it Let $\Gamma$ be any finitely generated
  group.  If $\Gamma$ is quasi-isometric to $\PZp$, then there is a
  short exact sequence
$$1 \rightarrow N \rightarrow \Gamma \rightarrow \Lambda \rightarrow
1$$
where $N$ is a finite group and $\Lambda$ is abstractly commensurable
to \linebreak $\PZp$.} 

\medskip
\noindent
Two groups are {\it abstractly commensurable} if they have isomorphic finite
index subgroups.

\subsection{An outline of the proofs of Theorems A and B}
\label{sec:outline}

The group $\PZp$ is a nonuniform (i.e. non-cocompact) lattice in $G=\G$ under the diagonal
embedding sending a matrix $M$ to the pair $(M,M)$.  
The group $G$ acts on $\A \times T_p$, where $\A$ is the hyperbolic plane and $T_p$ is the Bruhat-Tits-Serre tree associated to
$PGL_2(\Qp)$.  
Let $f: \PZp \rightarrow \PZq$ be a quasi-isometry.  
The proofs of Theorems A and B both begin as follows.

\bigskip
\noindent
{\bf Step 1 (The geometric model).}  We construct a space $\Omega_p$ on
which $\PZp$ acts properly discontinuously and
cocompactly by isometries (hence by a result of Milnor and Svarc \cite{Mi}, $\PZp$ and $\Omega_p$ are quasi-isometric).  
The space $\Omega_p$ has a boundary consisting of horospheres of
$\HT$, each of which is a quasi-isometrically embedded copy of the
group $BS(1,p)$.  
The quasi-isometry $f: \PZp \rightarrow \PZq$ then induces a
quasi-isometry, also denoted $f$, from $\Omega_p$ to $\Omega_q$.

\bigskip
\noindent
{\bf Step 2 (The Boundary Detection Theorem).}
This theorem shows that for every horosphere boundary component of
$\partial \Omega_p$, there is a corresponding horosphere boundary
component of $\partial \Omega_q$ so that $f$
restricts to a quasi-isometry of horospheres.  
The proof of this theorem uses the Coarse Separation Theorem of \cite{FS} and
\cite{S2}, and the geometry of $\Omega_p$.

\medskip
\noindent
{\bf Remark.} We are now able to prove Theorem B.  
The initial quasi-isometry $f: \PZp \rightarrow \PZq$ induces a
quasi-isometry $f: \Omega_p \rightarrow \Omega_q$.  
From the Boundary Detection Theorem, we obtain a
quasi-isometry $\hat f: \sigma \rightarrow \tau$ by restriction, where
$\sigma$ and $\tau$ are horosphere boundary components of $\Omega_p$
and $\Omega_q$, respectively.  
By Step 1, the map $\hat f$ can be considered as a quasi-isometry of
Baumslag-Solitar groups, namely
$\hat f: BS(1,p) \rightarrow BS(1, q).$
From Theorem \ref{thm:BS:qi} \cite{FM} we conclude that $p = q$.

\bigskip
\noindent
The proof of Theorem A continues with the following
steps.  We are now considering a quasi-isometry $f: \Omega_p
\rightarrow \Omega_p$.

\bigskip
\noindent
{\bf Step 3 (The geometry of $\Omega_p$).}  
For any two horosphere boundary components $\sigma_1$ and $\sigma_2$
of $\Omega_p$, there is a unique line $l \subset T_p$ so that
$\sigma_1$ and $\sigma_2$ are a specified fixed distance apart in $\A
\times t$, for all vertices $t$ of $l$.  
This line is called the {\it closeness line} of $\sigma_1$ and
$\sigma_2$.  
The set of closeness lines of all horospheres is preserved under
quasi-isometry.  
This geometric result replaces the usual group-equivariance assumed in
Mostow-Prasad rigidity, and provides the structure necessary for Step 4.

\bigskip
\noindent
{\bf Step 4 ($S$-Arithmetic Action
    Rigidity).}   
Here we prove an $2$ dimensional $S$-arithmetic version of the Action Rigidity Theorem
of R. Schwartz \cite{S}.  
The Action Rigidity Theorem concludes that the map induced by $f$ on
the set of horospheres of $\partial \Omega_p$ (which is indexed by $\Q
\cup \{ \infty \}$) is given by an affine map of $\Q \cup \{ \infty \}$.


\bigskip
\noindent
{\bf Step 5 (Conclusion of the proof).}  
From Step 4, we are able to choose a specific commensurator $g$ of $\PZp$ so that the composite map $f
\circ g$ is a bounded distance from the identity map.  This finishes
the proof of Theorem A.

\section{Preliminary material}
\label{sec:prelim}

\subsection{Quasi-isometries}
\label{sec:qi}

{\bf Definition.}  Let $K \geq 1$ and $C \geq 0$.  
A {\it $(K,C)$-quasi-isometry} between metric
spaces $(X, d_X)$ and $(Y,d_Y)$ is a
map $f:X \rightarrow Y$ satisfying:

\medskip
\noindent
1. $\frac{1}{K} d_X(x_1,x_2) - C \leq d_Y(f(x_1),f(x_2)) \leq K
d_X(x_1,x_2) + C$ for all $x_1,x_2 \in X$.

\medskip
\noindent
2. For some constant $C'$, the $C'$ neighborhood of $f(X)$ is all of
$Y$.

\bigskip

We often omit the constants $K$ and $C$ and simply refer to $f$ as a
quasi-isometry.  
A quasi-isometry $f$ can always be changed by a bounded amount using
the standard ``connect-the-dots'' procedure so that it is continuous.
(See, e.g. \cite{FS}.)
A quasi-isometry also has a coarse inverse, i.e. there is a quasi-isometry $g:Y
\rightarrow X$ so that $f \circ g$ and $g \circ f$ are a bounded
distance from the appropriate identity map in the sup norm.  
A map satisfying $1.$ but not $2.$ in the definition above is called a
{\em quasi-isometric embedding}.

We define the {\em quasi-isometry group} of a space $X$, denoted $QI(X)$, to
be the set of all self quasi-isometries of $X$, modulo those a bounded
distance from the identity in the sup norm, under composition of quasi-isometries.  
Inverses exist in $QI(X)$ since every quasi-isometry has a coarse
inverse.  
A quasi-isometry between two metric spaces $X$ and $Y$ induces an isomorphism
between $QI(X)$ and $QI(Y)$.

\subsection{$PSL_2$ as a algebraic group}
\label{sec:algebraic}

We will use the following model of $PSL_2$ as an algebraic group.
Consider the map $Ad: SL_2(\Cp) \rightarrow GL(sl_2)$, where we view
the Lie algebra $sl_2$ as a vector space.   
Let $G' = Ad(SL_2(\Cp))$.  
Then $G'$ is a model for $PSL_2$ as an algebraic group, since the center of $SL_2(\Cp)$ vanishes under the map $Ad$.  
By $PSL_2(\Q)$ we mean the $\Q-$points of $G'$, denoted $G'_{\Q}$.

\subsection{The geometry of $BS(1,n)$}
\label{sec:BS:geom}

The Baumslag-Solitar group $BS(1,n)$ acts properly discontinuously and cocompactly by
isometries on a
metric $2-$complex $X_n$ defined explicitly in \cite{FM}.  
This complex $X_n$ is topologically $T_n \times R$, where
$T_n$ is a regular $(n+1)$-valent tree, directed so that each vertex has
$1$ incoming edge and $n$ outgoing edges.  (Figure $1$.)

\begin{figure}
\begin{center}~
\epsfig{file=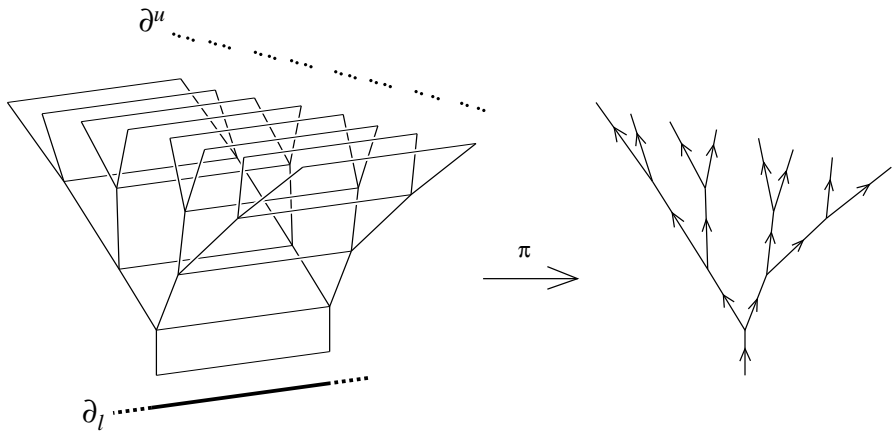,height=2.3in}
\caption[The $BS(1,n)$ complex]{The complex $X_n$ associated to $BS(1,n)$ where the map $\pi$
denotes projection onto the tree factor.  The upper and lower boundaries
are also marked.}
\end{center}
\end{figure}

A {\em height function} on $T_n$ is a continuous function $h:T_n \rightarrow \R$ which
maps each oriented edge of $T_n$ homeomorphically onto an oriented
interval of a given length $d$.  
A vertex of $T_n$ whose height under $h$ is $kd$ is said to have
{\em combinatorial height} $k$.  
Fix a basepoint for $T_n$ with height $0$.  This determines a height
function on $T_n$. 
Let $\pi: T_n \times \R \rightarrow T_n$ denote projection.  
Then $h \circ \pi$ is a height function on $X_n$.

A {\em proper line} in $T_n$ is the image of a proper embedding $\R
\rightarrow T_n$.  
A coherently oriented proper line is one on which the height function
is strictly monotone.  
We will use the term ``line'' to mean a proper line in $T_n$.   
The metric on $X_n$ is defined so that for each infinite,
coherently oriented line $l \subset T$, the plane $l \times \R$ is isometric to
a hyperbolic plane.

When studying the geometry of $X_n$, there are two ``boundaries'' of
the complex which play an important
role.  (Figure $1$.)  
The {\em lower boundary}, denoted $\partial_lX_n$, is homeomorphic to $\R$ and
is the common lower boundary of all hyperbolic planes in $X_n$.

The {\em upper boundary}, denoted $\partial^u X_n$, is defined to be the space of hyperbolic planes in
$X_n$, with the following metric.  
If $Q_1$ and $Q_2$ are hyperbolic planes in $X_n$ which agree below
combinatorial height $k$, define the distance between them to be
$n^{-k}$.  
With this
metric, $\partial^uX_n$ is isometric to the set of $n-$adic rational
numbers, $\Q_n$, with the metric
defined by the $n-$adic absolute value.

B. Farb and L. Mosher obtain the following quasi-isometric rigidity
results for $BS(1,n)$. \cite{FM} 

\begin{theorem}[\cite{FM}]
\label{thm:BS:qi}
For integers $m,n \geq 2$, the groups $BS(1,m)$ and $BS(1,n)$ are
quasi-isometric if and only if they are commensurable.  This happens
if and only if there exist integers $r,j,k>0$ such that $m=r^j$ and $n
= r^k$.
\end{theorem}

\begin{theorem}[\cite{FM}]
\label{thm:BS:2}
The quasi-isometry group of $BS(1,n)$ is given by the following isomorphism:
$$QI(BS(1,n)) \cong Bilip(\R) \times Bilip(\Q_n).$$
\end{theorem}

A quasi-isometry $f \in QI(BS(1,n))$ induces bilipschitz maps $f^u$ and $f_l$ on
the upper and lower boundaries of $X_n$, respectively \cite{FM}.    
From Theorem \ref{thm:BS:2} we see that the map $QI(BS(1,n)) \rightarrow Bilip(\R)
\times Bilip(\Q_p)$ given by $f \rightarrow (f_l, f^u)$ is an
isomorphism.  
It is perhaps surprising that $\PZp$ should have such a small
quasi-isometry group while $BS(1, p)$ has such a large quasi-isometry group.

\section{The geometry of $\PZp$}
\label{sec:geom}

\subsection{The action of $\G$ on $\HT$}
\label{sec:action}

We will consider $\PZp \subset \G$ as the image of the diagonal map $\eta: \PZp \rightarrow \PR \times
\PQp$ given by $\eta(M) = (M,M)$.  
Viewed in this way, $\PZp$ is a lattice in the group $\G$, for any
prime $p$.

We now define the Bruhat-Tits tree $T_p$ associated to $PGL_2(\Qp)$.  
We consider a tree $T$ to be a set of vertices $Vert(T)$ together with
a set of adjacency relations among the vertices.  
Let $Vert(T_p)$ be the set of equivalence classes of $\Zp$-lattices in
$\Qp \times \Qp$.  
Two lattices $L_1$ and $L_2$ are equivalent if $L_2 = \alpha L_1$,
where $\alpha \in \Qp - \{ 0 \}$.    
Two vertices $[L_1]$ and $[L_2]$ are adjacent if there exist
representatives $L_1 \in [L_1]$ and $L_2 \in [L_2]$ with $L_1 \subset
L_2$ and $[L_2 : L_1] = p$.  
An example of two adjacent vertices is $[\Zp \times \Zp]$ and $[p \Zp
\times \Zp]$.  
The tree $T_p$ is a regular $(p+1)$-valent tree. 
We will fix $[L_0] = [\Zp \times \Zp]$ as the basepoint of $T_p$ as
well as a height function $h$ giving $[L_0]$ height $0$.  
An element $g \in PGL_2(\Qp)$ acts on $[L] \in Vert(T_p)$  by matrix
multiplication on the basis vectors of a representative lattice in the
equivalence class $[L]$.  
Note that $PSL_2(\Qp)$ also acts on $T_p$ in this way.

Let $\A$ denote $2-$dimensional hyperbolic space in the upper half plane model, i.e. $\A = \{
(x,y) | x \in \R, y>0 \}$ with the metric $\frac{dx^2 + dy^2}{y^2}$.  
We define the action of an element $(g_1,g_2) \in \G$ on a
point $(x,[L]) \in \HT$.  
The element $g_1 \in \PR$ acts on $x \in \A$ by fractional linear
transformations.  
The element $g_2 \in \PQp$ acts on $[L] \in Vert(T_p)$ by matrix
multiplication on the basis vectors of a representative lattice in the
equivalence class $[L]$.  
When we are considering the action of $\PZp$ on $\HT$, we have $g_1
= g_2$.  
 We will use only one coordinate to represent the elements of $\PZp
\subset \G$.  
So $g \in \PZp$ would correspond to $(g,g) \in \G$.  
Hence we can refer to $g \in \PZp$ acting on either $\A$ or $T_p$ or
$\HT$ in the appropriate manner.  

\subsection{Constructing the space $\Omega_p$}
\label{sec:omegap}

We want to construct a space $\Omega_p \subset \HT$ on which
$\PZp$ acts properly discontinuously and cocompactly by isometries.  
The Milnor-Svarc criterion states that if a finitely generated group $\Gamma$ acts
properly discontinuously and cocompactly by isometries on a space $X$,
then $\Gamma$ is quasi-isometric to $X$.  
We then refer to the geometry of $X$ as the large scale geometry of
$\Gamma$.  
This additional geometric information associated to $\Gamma$ is
often useful in determining rigidity properties of $\Gamma$.

Although $\PZp$ acts by isometries on $\HT$, it does not act
cocompactly, because the fundamental domain for the action of $\PZp$
on a fixed $\A$ is the same
as the fundamental domain for the action of $PSL_2(\Z)$ on $\A$, which is
unbounded in one direction.

Let $w$ be the segment of the horocircle based at $\infty$ of height
$h_0$ in this fundamental domain (in the upper half space model of
hyperbolic space), for $h_0$
sufficiently large.  Fix $H>1.$  
Lift the segment $w$ to $\A \times [L_0]$ to obtain a horocyclic
segment at height $H$ whose orbit under $PSL_2(\Z)
= Stab_{\PZp} ([L_0])$ is a disjoint collection
of horocircles of $\A$, centered at $\Q \cup \{ \infty \} \subset \partial_{\infty} \A$.  
The orbit of the lift of $w$ under the entire group $\PZp$ gives
a $\PZp$-equivariant collection of horocircles in $\HT$, based at $\Q
\cup \{ \infty \}$ in each copy of $\A \subset \HT$.

We now define a {\em horosphere} of $\HT$ based at $\alpha \in \R \cup \{
\infty \}$ to be the collection of horocircles, one in each $\A \times
[L]$ for every $[L] \in Vert(T_p)$, all based at $\alpha$.  
According to the above construction, there is exactly one such
horocircle in each $\A \times [L]$.  
We will denote this horosphere of $\HT$ by $\sa$.

In order to have a ``connected'' picture of a horosphere, we can put
an edge $e$ between any two adjacent vertices of $T_p$ and extend the
horosphere linearly in $\A \times e$.  
So we can think of a horosphere as a hollow tube, or in the case of
$\si$ as a flat sheet, whose image under the projection $\pi: \HT
\rightarrow T_p$ is all of $T_p$. (Figures $2$ and $3$.)

We define the space $\Omega_p$, where $\PZp$ acts properly
discontinuously and cocompactly by isometries, to be
$\HT$ with the interiors of all the horospheres removed.  The interior
of a horosphere is the union of the interiors of the component
horocircles.

\begin{figure}
\begin{center}~
\epsfig{file=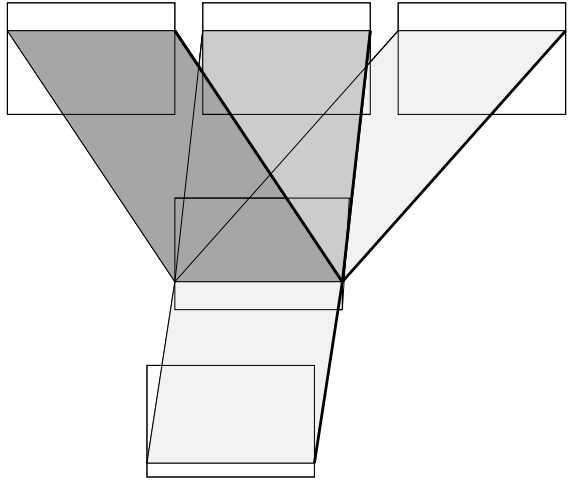,height=2.5in}
\caption[The horosphere based at $\infty$ in $\Omega_p$]{This is a
  piece of the horosphere $\si$ in $\A \times T_3$, the space associated to $PSL_2(\Z[\frac{1}{3}])$.  
Notice how the height of the horosphere increases in
  each successive copy of $\A$.  
Viewing the bold black lines as part of the tree $T_3$ helps one to see
that this horosphere is topologically $T_3 \times \R$.}
\end{center}
\end{figure}

\subsection{The metric}
\label{sec:metric}

The following theorem allows us to use the product metric $d_{\A}
\times d_T$ on $\Omega_p$.  
Although it is stated for semisimple Lie groups, it is proven in the
more general context of $S$-arithmetic groups.  

\begin{theorem} \cite{LMR}
If $G$ is a semisimple Lie group of rank at least $2$ and $\Gamma$ is
an irreducible lattice in $G$ then, $d_R$ restricted to $\Gamma$ is
Lipschitz equivalent to $d_W$, where $d_W$ is the word metric on
$\Gamma$ and $d_R$ is the left invariant Riemannian metric on
$\Gamma$.  
\end{theorem}

By construction, $\PZp$ acts properly discontinuously and cocompactly
by isometries on $\Omega_p$.

\subsection{The closeness line}
\label{sec:closeness:line}

We are now interested in the packing of the horospheres of $\HT$, i.e. how the
horosphere ``tubes'' fit together. 
We consider two horospheres, without loss of generality $\so$ and
$\si$, and the distance between them when restricted to $\A \times
[L]$, for any $[L] \in Vert(T_p)$.  
We use the notation $\sigma|_{[L]}$ for the horosphere $\sigma \subset
\partial \Omega_p$ restricted to $\A \times [L]$.  
We will show that there is a unique line $l \subset T_p$ so that $\so$
and $\si$ remain a constant distance apart when restricted to $\A
\times [L]$, for all vertices $[L]$ lying on $l$.  
In addition, we will show that all other lines $l' \subset T_p$ have the property that
$d_{\A}(\so|_{[L]}, \si|_{[L]})$ increases without bound as the vertices $[L]$
lying on $l'$ increase in height. 
We call $l$ the {\em closeness line} of $\si$ and $\so$.  
By symmetry, any two horospheres $\sa$ and $\sigma_{\beta}$ have a
unique closeness line.

We will use the notation $\sigma|_{[L]}$ for the horosphere $\sigma$ of
$\partial \Omega_p$ restricted to $\A \times [L]$, for $[L] \in
Vert(T_p)$.  
Note that a matrix $g$ which translates $\sa|_{[L_1]}$ to
$\sa|_{[L_2]}$ will lie in $Stab_{\PZp}(\alpha)$ and satisfy $[gL_1] =
[L_2]$.    
However, to move between the horocircles $\sigma_{\beta}|_{[L_1]}$
and $\sigma_{\beta}|_{[L_2]}$, for $\beta \neq \alpha$, we will need to use a different
element of $\PZp$, i.e. $h \in Stab_{\PZp}(1)$ such that $[hL_1] =
[L_2]$.

Products of the matrices $ A = \left(
\begin{array}{cc} p & 0 \\ 0 & \frac{1}{p} \end{array} \right)$ and  $
    B = \left(
\begin{array}{cc} 1 & 1 \\ 0 & 1 \end{array} \right)$ give the
$PSL_2(\Qp)$ action on $T_p$. (See, e.g. \cite{Se}.)
Since 
$$A \in Stab_{\PZp}(0) \cap Stab_{\PZp}( \infty),$$ letting $A^i$
act on $\A \times [L_0]$ moves the horocircle $\so|_{[L_0]}$ (resp. $\si|_{[L_0]}$) to the horocircle  $\so|_{[A^iL_0]}$  (resp. $\si|_{[A^iL_0]}$).  
Moreover, $A \in PSL_2(\R) = Isom ^+(\A)$, so we have 
$$d_{\A}(\so|_{ [L_0]}, \si|_{[L_0]}) = d_{\A}(\so|_{[A^iL_0]}, \si|_{[A^iL_0]}) = 2 \log H
$$
where $H$ was chosen in \S \ref{sec:omegap}.  
Let $l$ be the line in $T_p$ which is the orbit of $[L_0]$ under the
cyclic group generated by the matrix $A$.  
We call $l$ the {\em diagonal line}.  

Recall that we usually need two different matrices to move two
different horocircles in $\A \times [L_0]$ to their corresponding
horocircles in $\A \times [L]$ for $[L] \in Vert(T_p)$; here we need
only one matrix since the matrix $A$ lies in the intersection of the
stabilizers of $0$ and $\infty$.  

We will view a ray of $T_p$ based at $[L_0]$ as an infinite sequence of
products $\{ \Pi_{i=1}^N C_i \}_{N \in \N}$ where either $C_i = A$ or
$C_i = B^jA$.  
The diagonal ray (i.e. the part of the diagonal line beginning at
$[L_0]$ and moving upwards in height) is  described by $\{ \Pi_{i=1}^N
A \}_{N \in \N}$, and we know from the previous paragraph that 
$d_{\A}( \so|_{[\Pi_{i=1}^N A]}, \si|_{[\Pi_{i=1}^N A}]) = 2 \log H$ for
all $N \in \N$.  

Now suppose that $r \subset T_p$ is the ray based at $[L_0]$ given by
the sequence $\{ \Pi_{i=1}^N C_i \}_{N \in \N}$ where $C_1 = B^jA$ for
$j \in \{ 0,1, \cdots p-1 \}$.
For any $N_0 \in N$, the product $\Pi_{i=1}^N C_i$ has the form
$(\Pi_{i=1}^{T} A)( \Pi_{i=1}^S C_i)$ where $C_i$ is as above.  
Such a product is a matrix of the form $M = \left( \begin{array}{cc} p^{S+T} & \frac{s}{p^{S-T}} \\ 0 &
\frac{1}{p^{S+T}} \end{array} \right)$, where $S$ or $T$ may be $0$. 
We may assume that $S>T$, meaning that in the tree
factor, we are
considering vertices sufficiently far away from the diagonal line
$l$.

Let $\beta = \frac{-s}{p^{2S}}$, where $s$ and $S$ are as in the matrix $M$.  
We can find a matrix of the form $N = \left( \begin{array}{cc} a & -s \\ b &
p^{2S} \end{array} \right) \in PSL_2(\Z)$ where $s, S$ are as in $M$.  
Since $N \in PSL_2(\Z) = Stab_{\PZp}([L_0])$ and $N \cdot 0 =
\beta$, the matrix $N$ maps $\so|_{[L_0]}$ to $\sigma_{
  \beta}|_{[L_0]}$.  
A computation shows that the height of $\sigma_{
  \beta}|_{[L_0]}$ is $\frac{1}{p^{2(S+T)}H}$.  
Moreover, $M \cdot \beta = 0$, so $M \cdot (\sigma_{\beta}|_{[L_0]}) =  \so|_{[ML_0]}$.  
Another computation shows that $\so|_{[ML_0]}$ has height
$\frac{1}{p^{2m}H}$ in $\A$.  
The height of $\si|_{[ML_0]}$ is $p^{2(S+T)}H$.  
Hence 
$$
d_{\A}(\so|_{[ML_0]}, \si|_{[ML_0]}) = \log (p^{2S}H^2).$$  
When $S=0$, the vertex $\Pi_1^T A$ lies on the diagonal line, and the
above formula gives a constant distance of $2 \log H$ between $\si$ and
$\so$ in $A \times l$.  
When $S \neq 0$, we see that $S$ increases by a factor of $2$ for each
unit of height in $T_p$; hence the distance between $\si$ and $\so$
increases by a factor of $p^4$.

So we see that for any line $l' \subset T$ which is not the diagonal
line, \linebreak 
$d_{\A}(\so|_{[L]}, \si|_{[L]})$ increases without bound as we
consider vertices $[L]$ of $l'$ of increasing height.  
So the diagonal line $l$ is the closeness line of $\so$ and
$\si$.

It is clear by construction that the closeness line of $\sa$ and
$\sigma_{\beta}$ and the closeness line of $\sa$ and $\sigma_{\gamma}$
are distinct.

\subsection{The geometry of the horospheres}
\label{sec:BS}

Consider the matrices $A = \left( \begin{array}{cc} p & 0 \\ 0 & \frac{1}{p}
\end{array} \right)$ and $B = \left( \begin{array}{cc} 1 & 1 \\ 0 & 1
\end{array} \right)$.  
Then $A,B \in \PZp \subset \PQp$ and   
$$ABA^{-1} = \left( \begin{array}{cc} 1 & p^2 \\ 0 & 1 \end{array}\right) =
\left( \begin{array}{cc} 1 & 1 \\ 0 & 1
\end{array} \right)^{p^2} = B^{p^2}.$$
The map $\Phi$ sending the generator $a$ of $BS(1,p^2)$ to the matrix $A$ and
the generator $b$ to $B$ is a homomorphism of $BS(1,p^2)$ into the group $\PQp$.  
Products of the matrices $A$ and $B$ form the elements of $PGL_2(\Qp)$
needed to move between vertices of $T_p$.   

If we take the orbit of the segment $w$ from $Hi$ to $1 + Hi$ in $\A
\times [L_0]$ under the group $\Phi(BS(1,p^2)) \subset \PQp$, we
obtain the horosphere $\si$.  
The width of a horostrip in $\si$ between vertices whose combinatorial
heights differ by $1$, in the metric
on $\HT$, is $1 + 2 \log p$, while in $X_{p^2}$ it is $2 \log p$.
Since these distances are comparable, the horosphere $\si$ (and hence
any horosphere $\sa$) is quasi-isometric to the complex $X_{p^2}$
associated to $BS(1,p^2)$.  
Since $BS(1,p^2)$ and $BS(1,p)$ are commensurable, hence
quasi-isometric, we may assume
that the horospheres are quasi-isometrically embedded copies of the
complex $X_p$ associated to $BS(1,p)$.

\subsection{Another view of the closeness lines}
\label{sec:geom:closeness}

The following discussion of closeness lines provides some geometric
intuition for understanding these lines but is not necessary for
the proofs of Theorems A, B and C.

We can consider the collection of closeness lines of all
horospheres of $\partial \Omega_p$ with the horosphere $\si$.  This is
a collection of distinct lines in the tree $T_p$.  
These closeness lines can also be viewed as lying in the horosphere
$\si$, which is thought of as the complex $X_p$ associated to $BS(1,p)$.  
In this case, the closeness line of $\si$ and $\sa$ can be viewed as
the set of points in $\si$ closest to $\sa$ which project under $\pi:
\HT \rightarrow T_p$ to the closeness line of $\si$ and $\sa$ as
described in \S \ref{sec:closeness:line}.  (Figure $4$.)  
Then each hyperbolic plane $\A \subset X_p = \si$ intersects these closeness lines, with at most one line lying completely in
that plane.  (Figure $5$.)

\begin{figure}
\begin{center}~
\epsfig{file=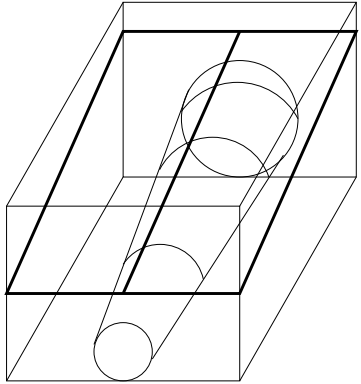,height=2in}
\caption[Closest points in $\si$]{The two horospheres shown are $\sa$ and
  $\si$, for $\alpha \neq \infty$.  The dark line in $\sigma_{\infty}$ represents
  the closest points in $\sigma_{\infty}$ to $\sa$.  We can view
  this line as the closeness line of $\sigma_{\infty}$ and
  $\sa$ drawn in $\sigma_{\infty}$.}
\end{center}
\end{figure}

\begin{figure}
\begin{center}~
\epsfig{file=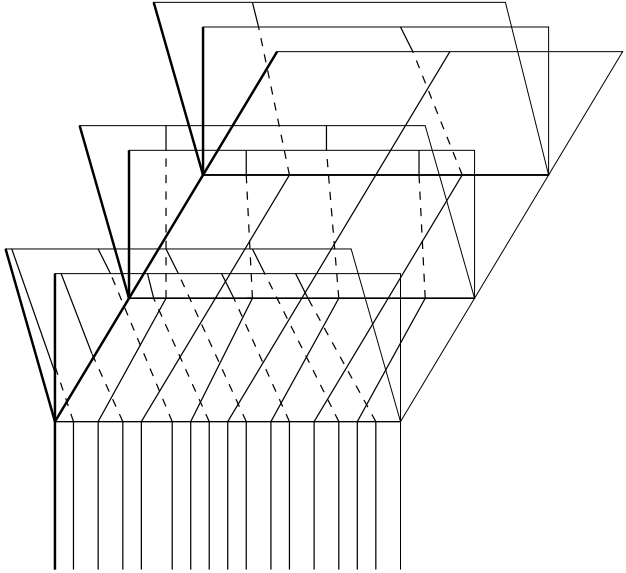,height=3in}
\caption[Closeness lines in $\si$]{A subset of the collection of closeness lines of all
  horospheres with $\sigma_{\infty}$, drawn in the complex $X_p$ of $\sigma_{\infty}$.}
\end{center}
\end{figure}

\section{The Boundary Detection Theorem}
\label{sec:horospheres}

The first goal of this section is to prove Theorem \ref{thm:BDT}
(Boundary Detection Theorem),
which states that a quasi-isometry $f: \Omega_p \rightarrow \Omega_q$
maps a horosphere boundary component of $\Omega_p$ to within a bounded
Hausdorff distance of a horosphere boundary component of $\Omega_q$.  
Two subsets $X$ and $Y$ of a metric space $W$ have bounded Hausdorff
distance if there exists an $\epsilon > 0$ so that $X \subset
Nbhd_{\epsilon} (Y)$ and $Y \subset
Nbhd_{\epsilon} (X)$. 
The infimum of such constants $\epsilon$ is called the {\em Hausdorff distance} between $X$
and $Y$.  

The second goal is to use Theorem \ref{thm:BDT}, combined
with Theorem \ref{thm:BS:qi}, to prove Theorem B, i.e. that $\PZp$ and $\PZq$ are quasi-isometric
if and only if $p = q$.  

We use the notation $Nbhd_r(S)$ to be the $r$-neighborhood in $\Omega_p$ of
a subset $S$ of $ \Omega_p$.  
We say that a subset $S$ of a metric space $X$ has the
{\em strong separation property} in $X$ if there is  a fixed $r>0$
with the following property.  
For every $k > 0$, there are at least two connected components of $X -
Nbhd_r(S)$ which contain metric balls of radius $k$.  
We will say that $S$ {\em separates $X$} if $S$ has the strong separation
property in $X$.  The constant $r$ is called the {\em separation constant}.

A metric space $X$ is called {\em uniformly contractible} if there is a function
$\alpha: \R^+ \rightarrow \R^+$ with the following property.  
If $f : \Delta \rightarrow X$ is a map of a finite simplicial complex,
and $f(\Delta) \subset B_r$, where $B_r \subset X$ is a metric
$r$-ball, then $f(\Delta)$ is contractible in $B_{\alpha(r)}$, where
$\alpha$ is independent of $dim(\Delta)$.  
Any contractible space admitting a cocompact group of isometries is
uniformly contractible.

We will need the following coarse topology results which we state as
special cases of Theorem $5.2$ and Corollary $5.3$ of \cite{FS}.  
We are using the bounded geometry metric $d$ on $\R^3$ described below which comes from choosing a proper embedding of the tree $T_p
\rightarrow \R^2$.   This allows us to consider $T_p \times \R$ as a
subset of $(\R^3,d)$.  
In applying the results of \cite{FS} we use the fact that $\R^3$ in
this bounded geometry metric is uniformly
contractible and contains spheres of arbitrarily large radius
(the ``expanding spheres'' condition of \cite{FS}).

\begin{theorem}[Coarse Separation \cite{FS}]
\label{thm:CST}
Suppose $\phi: (\R^3,d) \rightarrow Y$ is a $(K,C)-$quasi-isometric embedding of
$(\R^3,d)$ into a uniformly contractible Riemannian manifold $Y$
diffeomorphic to $\R^4$.  Then $\phi(\R^3)$ separates $Y$, where the
separation constant depends on $(K,C)$.  
\end{theorem}

\begin{corollary}[Packing Theorem \cite{FS}]
Suppose $\phi: (\R^3,d) \rightarrow (\R^3,d)$ is a $(K,C)-$quasi-isometric
embedding.  Then $\phi$ is a $(K',C')-\linebreak $quasi-isometry, for some
constants $(K',C')$ depending on $(K,C)$.
\end{corollary}

\subsection{Separation}
\label{sec:separation}

Let $\sigma$ be any horosphere boundary component of $\Omega_p$ and
$f: \Omega_p \rightarrow \Omega_q$ a quasi-isometry.  In this
section we show that the image $f(\sigma)$ separates $\HTq$.  
To prove this, we extend the quasi-isometric embedding $f|_{\sigma}:
\sigma \rightarrow \HTq$ to a quasi-isometric embedding
$\widehat{f}: (\R^3,d) \rightarrow \A \times \R^2$ to which we can apply
Theorem \ref{thm:CST} (Coarse Separation). 
As above, $d$ is the bounded geometry metric on $\R^3$ which comes from choosing a proper embedding of the tree $T_p
\rightarrow \R^2$.

Consider $\sigma$ as the complex $X_{p}$ associated to $BS(1,p)$ and choose a homeomorphism
$\beta: X_{p} \rightarrow T_p \times \R$.  
Let $\alpha_p: T_p \rightarrow \R^2$ be any proper embedding.  
Then we can consider $X_{p}$ as a subset of $\R^3$ via the map
$(\alpha \times Id) \circ \beta$.  
It is shown in \cite{FM} that $\R^3$ can be given a bounded geometry
metric $d$ for which this map is an isometric embedding, as follows.  
The boundary
of each connected component $C$ of $\R^3 - X_{p}$ is topologically a
plane.  
We use two coordinates $(t,r)$ on this plane $\partial C$, where $t \in T_p$ and $r
\in \R$.  
Choose a homeomorphism which identifies $ \partial C \cup C$ with
$\partial C \times [0, \infty)$.  
Then a point in $C$ has three coordinates: $(t,r,s)$ where $t$ and $r$
as above give a point in $\partial C$ and $s \in [0, \infty)$.  
We use the product metric on each component.   

The quasi-isometric embedding $f|_{\sigma}$ has image in $\HTq$.  
Analogous to the above situation, we choose a proper embedding
$\alpha_q: T_q \rightarrow \R^2$ and view $\HTq$ as a subset of $\A
\times \R^2$ via the map $Id \times \alpha$.  
As above, we obtain a metric $d'$ on $\A \times \R^2$ for which this
map is an isometric embedding.   
To apply Theorem \ref{thm:CST} we use the fact that $(\A \times
\R^2,d')$ is diffeomorphic to $\R^4$, uniformly contractible and contains
arbitrary large metric balls.  

We can now extend the map $f|_{\sigma}$ to a quasi-isometric embedding $$\widehat f : (\R^3,d)
\rightarrow (\A \times \R^2,d')$$ 
by $$\widehat f(t,r,s) = (f(t,r),s).$$
Note that $f(t,r)$ is a point in $\HTq$, and so provides two
coordinates.  Also, $\widehat f|_{\sigma} = f$.   

\begin{proposition}[Separation]
\label{prop:sep}
Let $f: \Omega_p \rightarrow \Omega_q \subset \HTq$ be a $(K,C)-$quasi-isometry, 
and let $\sigma \subset \partial \Omega_p$ be any horosphere boundary component.  Then
$f(\sigma)$ separates $\HTq$, where the separation constant depends on
$(K,C)$.  
\end{proposition}

\begin{proof}
Extend $f|_{\sigma}$ as above to a quasi-isometric embedding $$\hat f: (\R^3,d)
\rightarrow (\A \times \R^2,d').$$  
It is understood that a path ``avoiding'' $f(\sigma)$ or $\widehat f(
\R^3)$ stays outside the neighborhood $Nbhd_r(f(\sigma))$ or $Nbhd_{r'}(\widehat f(
\R^3))$ for a constant $r$ or $r'$.  
Also note that $f(\sigma) = \widehat f(\sigma)$. 

The map $\widehat f: (\R^3,d) \rightarrow (\A \times \R^2,d')$ satisfies the
conditions of Theorem \ref{thm:CST}, so $\widehat f(\R^3)$
separates $\A \times \R^2$. 

Suppose that  $\hat f( \R^3)$ separates $\A \times \R^2$ and $\widehat f(\sigma)$ does not separate $\HTq$.  
This means that any two points in $ (\HTq) - \widehat f(\sigma)$ can be joined by a path avoiding $\widehat f(\sigma)$.    
We will show that if $\widehat f(\sigma)$ does not separate $\HTq$, any two points in $(\A \times \R^2)  - \widehat{
  f} ( \R^3 )$ can be connected by a path avoiding $\widehat{ f}
( \R^3 )$, contradicting the fact that $\widehat f (\R^3)$ separates
$\A \times \R^2$.  
Let $x_1$ and $x_2$ be any two points in $( \A \times \R^2) - \widehat{ f} ( \R^3 )$.  
Each $x_i, i=1,2$, has coordinates $(\alpha_i, t_i,s_i)$ as above.  
Since $x_i \notin \widehat f(\R^3)$, each line $\{ (\alpha_i, t_i, s )|
s \in \R \}$ is not contained in $\widehat f(\R^3)$ by
construction.  
When $s=0$, each line gives a point in $\HTq$ not contained in $\widehat f(\sigma)$ .  
Call these points $\beta_i$.
Since $\widehat f(\sigma)$ does not separate $\HTq$, we can connect $\beta_1$ to
$\beta_2$ by a path $\gamma$ lying in $\HTq$ which avoids
$\widehat f(\sigma)$ and hence $\widehat f(\R^3)$.  
So $x_1$ and $x_2$ are connected by the path 
$$( \alpha_1,t_1, [0,s_2])
* \gamma * ( \alpha_2, t_2, [0,s_1])^{-1}$$
 which avoids $\widehat{ f} ( \R^3 )$.  
Thus,  if $\widehat{ f} ( \R^3 )$ separates $\R^2  \times \hyp^2$ then
$f(\sigma) = \widehat f(\sigma)$ separates $\HTq$. \end{proof}

We now prove a lemma which shows that the space $\Omega_p$ with the
neighborhood of any horosphere removed is path connected.  
This lemma is needed in the proof of Theorem \ref{thm:BDT}. 

\begin{lemma}
\label{lemma:unobtrusive:boundaries}

Let $r$ be any positive real number.  For any horosphere boundary
component $\sigma$ of $\Omega_p$, the space $\Omega_p - Nbhd_r(\sigma)$ is path
connected.  
\end{lemma}

\begin{proof}
We will show that $\Omega_p - Nbhd_r(\sigma)$ is path connected for any $\sigma
\subset \partial \Omega$ and for any $r$.  
For convenience we may assume that $\sigma = \si$.  

For any horosphere $\tau \subset \partial \Omega_p$, consider the set
of vertices of $T_p$ where $\tau$ intersects $Nbhd_r(\si)$:
$$S_{r, \tau} = \{ [L] \in Vert(T_p) | \ \tau|_{[L]} \cap Nbhd_r(\si)
  \neq \emptyset \}.$$
Since each $\tau$ has a closeness line with $\si$, for large enough
$r$ the set $S_{r,
  \tau}$ is the neighborhood in $T_p$ of a line $l_{\tau} \subset T_p$.  
There is one line $l_{\tau}$ for each horosphere $\tau \neq \si$.  
Since the horosphere boundary components of $\Omega_p$ are indexed by
rationals and the set of lines of $T_p$ is indexed by elements of
$\Q_p$, there are lines $l \subset T_p$ such that no horosphere intersects
$Nbhd_r(\si)$ over every point of $l$.  Let $l'$ be such a
line.  

Choose two points $x$ and $y$ in $\Omega_p - Nbhd_r(\sigma)$.  If $x$ and
$y$ lie in $\A \times l'$, then there is a path connecting them.  
If this is not the case, let $\alpha$ be any path from $x$ to a point
$x' \in \A \times l'$, and $\beta$ any path from $y$ to a point $y'
\in \A \times l'$.  Let $\gamma$ be a path in $\A \times l'$ connecting $x'$ and $y'$.  
The composite path $\beta^{-1} \circ \gamma \circ \alpha$ connects $x$
to $y$ and lies in $\Omega_p - Nbhd_r(\sigma)$, proving the theorem. 
\end{proof}

\subsection{Proof of the Boundary Detection Theorem}
\label{sec:pf:BDT}

We now state and prove Theorem \ref{thm:BDT}, which plays a role in
the proofs of Theorems A and B.

\begin{theorem}[Boundary Detection Theorem]
\label{thm:BDT}
Let $f: \Omega_p \rightarrow \Omega_q$ be a $(K,C)-$quasi-isometry.  
There exist constants $(K',C')$ depending on $(K,C)$ and the spaces $\Omega_p$ and $\Omega_q$, with the following
property.  
For every horosphere boundary component $\sigma$ of $\Omega_p$, there
is a horosphere boundary component $\tau$ of $\Omega_q$ so that
$f|_{\sigma}: \sigma \rightarrow \tau$ is a $(K',C')-$quasi-isometry.
\end{theorem}

The following lemmas form the two major components of the proof of
Theorem \ref{thm:BDT}. 
In these lemmas, $f: \Omega_p \rightarrow \Omega_q$ will be a $(K,C)-$quasi-isometry.  

\begin{lemma}
\label{lemma:BDT:1}
There exists a constant $\epsilon$ depending on $(K,C)$ and the
spaces $\Omega_p$ and $\Omega_q$, with the following property.  
For every horosphere boundary component $\sigma$ of $\Omega_p$, there
is a horosphere boundary component
$\tau$ of $\Omega_q$ so that $\tau \subset
Nbhd_{\epsilon}(f(\sigma))$.  
\end{lemma}

\begin{proof}
  
Consider a horosphere boundary component $\sigma$ of $\Omega_p$.  
We will show that any two points $x_1$ and $x_2$ of $\HTq$
which are outside $Nbhd_{\delta}(f(\sigma)) $ (for some constant
$\delta$) can be connected by a path avoiding $Nbhd_{\delta}(f(\sigma))
$.  
If the lemma is false, then any horosphere $\tau \subset \partial
\Omega_q$ contains points not within $\epsilon$ of $f(\sigma)$.
For the proper choice of $\epsilon$, this will contradict the fact that $f(\sigma)$ separates $\HTq$.  

Choose points $x_1, x_2 \notin
Nbhd_{\delta}(f(\sigma)) $.  
If $x_i \notin \Omega_q$, $i=1,2$, we can find a path connecting $x_i$ to some $x_i' \in
\Omega_q$ which avoids $Nbhd_{\delta}(f(\sigma))$.  
Hence we will assume that $x_i \in \Omega_q$.

Let $f^{-1}$ be a coarse inverse to $f$ and consider the points $y_i =
f^{-1}(x_i) \subset \Omega_p$ for $i = 1,2$.  
If $\delta $ is sufficiently large, these points lie outside
$Nbhd_{\delta_1}(\sigma)$, for some $\delta_1$.  
From Lemma \ref{lemma:unobtrusive:boundaries} we know that there is a path $\gamma$ between
$y_1$ and $y_2$ which avoids $Nbhd_{\delta_1}(\sigma)$.  
Then for some $\epsilon$, $f(\gamma)$ is a path in $\Omega_q$
connecting $x_1$ and $x_2$ but avoiding $Nbhd_{\epsilon}(f(\sigma))$.  
(Figure $7$.)  
If our initial $\delta$ was large enough, then $\epsilon > \delta$
and we can make this argument with $x_1 \in \tau, \ x_1 \notin
Nbhd_{\delta}(f(\sigma))$, contradicting the fact that $f(\sigma)$ separates $\HTq$.  \end{proof}

\begin{lemma}
\label{lemma:BDT:2}
There exists a constant $\epsilon'$ depending on $(K,C)$ and the
spaces $\Omega_p$ and $\Omega_q$, with the following property.  
For every horosphere boundary component $\sigma$ of $\Omega_p$, there is a horosphere boundary component $\tau$ of
$\Omega_q$ so that $f(\sigma) \subset Nbhd_{\epsilon'}(\tau)$.
\end{lemma}

\begin{proof}
Consider a horosphere boundary component $\sigma$ of $\Omega_p$.  
From Lemma \ref{lemma:BDT:1}, there is a horosphere boundary component
$\tau$ of $\Omega_q$ and a constant $\epsilon$ so that $\tau \subset Nbhd_{\epsilon}(f(\sigma))$.
Define a map $\psi: \tau \rightarrow \sigma$ by $\psi(y)  = x \in
\sigma$, where $x$ is any point so that $
f(x)$ is metrically closest to $y$. 
If there is more than one such point, choose randomly.   
From Lemma \ref{lemma:BDT:1}, we see that $\psi$ differs from the
coarse inverse $f^{-1}$ of $f$ by at most a constant.  
From Proposition \ref{prop:sep} we know that $\psi(\tau)$
separates $\HTq$.  
So for some constant $\delta'$, the horosphere $\sigma$ is contained
in $Nbhd_{\delta'}(\psi(\tau))$, i.e. every point of $\sigma$ is within
a constant $\delta'$ of some point $x \in \psi(\tau)$ which maps to a point $f(x)$
within $\epsilon$ of a point of $\tau$.  
By enlarging $\delta'$ to a constant $\epsilon'$, we see that $f(\sigma) \subset
Nbhd_{\epsilon'}(\tau)$.  Since all the horospheres are isometric, the constant $\epsilon'$ is
independent of the choice of $\sigma$.  
\end{proof}

\bigskip
\noindent
{\it Proof of the Boundary Detection Theorem.}  
Apply Lemma \ref{lemma:BDT:2} to both $f$ and $f^{-1}$.  
Compose $f$ with a nearest point projection to obtain a
quasi-isometric embedding $f:\sigma \rightarrow \tau$.  
As in \S \ref{sec:separation}, embed $\sigma$ and $\tau$
isometrically in $(\R^3,d)$ and extend $f$ to a quasi-isometric embedding $\hat f:
(\R^3,d) \rightarrow (\R^3,d)$. 
(Recall that $d$ was the bounded geometry metric on $\R^3$ described
in \S \ref{sec:separation}.)   
From the Packing Theorem, this map is a $(K',C')-$quasi-isometry,
where the pair $(K',C')$ depends on $(K,C)$.  
Then by the construction of $\hat f$, the map $f$ must
also be a $(K',C')-$quasi-isometry. 
\qed

It is now clear that closeness lines are preserved under
quasi-isometry.

\subsection{Proof of Theorem B}
\label{sec:qi:comm}

Theorem \ref{thm:BDT} allows us to prove Theorem B.

\bigskip
\noindent
{\it Proof of Theorem B.}  
If $\PZp$ and $\PZq$ are commensurable, then they are
automatically quasi-isometric.  

Let $f: \PZp \rightarrow \PZq$ be a $(K,C)-$quasi isometry.  
Construct the spaces $\Omega_p$ and $\Omega_q$
corresponding to $\PZp$ and $\PZq$, respectively.  
Then $f$ induces a quasi-isometry, also denoted $f$, from $\Omega_p$ to
$\Omega_q$.  
Theorem \ref{thm:BDT} allows us to restrict $f$ to a quasi-isometry $\hat f$ on
horospheres.  In \S \ref{sec:BS:geom}  we showed that a horosphere of
$\HT$ has the geometry of the group $BS(1,p)$.  
Hence $\hat f$ is a quasi-isometry from $BS(1,p)$ to $BS(1,q)$.  
According to Theorem \ref{thm:BS:qi}, we must have $p=q$ for
such a quasi-isometry to exist. 
\qed

\section{Theorems A and C}
\label{sec:main:theorem}

Every commensurator of $\PZp$ acts as a quasi-isometry of $\Omega_p$.  
To prove Theorem A, we must show that by
composing an element $f \in QI(\PZp)$ with a specific commensurator, we obtain a map which is a
bounded distance from the identity map.  
It is Theorem \ref{thm:action:rigidity} (stated below) that tells us which
commensurator to choose for this purpose.  
In the Appendix, we show that the
commensurator group of  
$\PZp$ in $\G$ is $PSL_2(\Q)$, where $PSL_2$ is viewed as an
algebraic group as in \S \ref{sec:algebraic}.  

In all that follows, we assume that $f: \Omega_p
\rightarrow \Omega_p$ is a $(K,C)-$quasi-isometry which has been
changed by a bounded amount using the ``connect
the dots'' procedure so that it is continuous. (See, e.g. \cite{FS}.)
Since $Comm(\PZp) = PSL_2(\Q)$
acts transitively on pairs of distinct points of $\R \cup \{ \infty \}$, we can assume that $f$ has been composed with a commensurator so that
$f(\si) = \si$ and $f(\so) = \so$.  
(Note that these horospheres are not necessarily fixed pointwise.)

\subsection{Action Rigidity}
\label{sec:background:ar}

We will now define a boundary of the space $\HT$.  
Let $t$ be the common endpoint of any two lines in $T_p$ and consider
$\partial_{\infty}(\A \times t) - \{ \infty \} \cong \R$.  
Recall from \S \ref{sec:BS:geom} that we can consider $\si \subset
\partial \Omega_p$ as a quasi-isometrically embedded copy of the
$2$-complex $X_p$ associated to $BS(1,p)$.  
So we can refer to the lower boundary $\partial_l(\si)$ of $\si$
(resp. the upper boundary $\partial^u(\si)$) as
the lower (resp. upper) boundary of $X_p$.  
The inclusion $i : \si \rightarrow \HT$ induces an identification
between $\partial_l(\si)$ and the copy of $\R$ described above.  
Since $f(\si) = \si$, the quasi-isometry $f$ induces bilipschitz maps $f_l$ and $f^u$ on
$\partial_l(\si)$ and $\partial^u(\si)$ of $\si$, respectively.  
The restriction of $f$ to $\R = \partial_{\infty}(\A \times t) - \{
\infty \}$ determines the permutation of the horospheres of $\HT$
under the map $f$.  
Hence, from the identification induced by the inclusion map $i$, the
map $f_l$ is exactly the map which determines the permutation of the
horospheres under $f$.

Let $\dq$ denote the lattice in $\R \times \Qp$ given by the diagonal $
\{ (a,a) | a \in \Q \} \subset \R \times \Qp$, and $\dz$ the sublattice $\{ (b,b) | b \in \F \}$.  
Clearly $\dz \subset \dq$.  We view of the first coordinate of
the pair $(b,b) \in \dq$ as $b \in \Q \subset \R$ denoting the basepoint of
a horosphere of $\HT$, and the second coordinate $b \in \Q \subset \Qp$ as the
point in $\Qp$ determined by the closeness line of $ \sigma_b$ and $\si$.  
Using the identification described above, we see that the map induced
by $f$ on $\Q \subset \R$ is exactly $f_l|_{\Q}$.  
We can view $\Q \subset \Q_p$ as a subset of $\partial^u(\si)$.  
Since $f(\si) = \si$, the map induced by $f$ on $\Q \subset \Q_p$ is
exactly $f^u|_{\Q}$.  
The maps $f^u|_{\Q}$ and $f_l|_{\Q}$ are identical because closeness
lines are preserved under quasi-isometry and $f(\si) = \si$.  
Hence $f$ induces a map of $\dq$ given by $(f_l|_{\Q}, f^u|_{\Q})$,
and we
will use a single coordinate for points of $\dq$.  
Let $\phi$ denote
the common restriction of $f_l$ and $f^u$ to $\Q$.  
Then $\phi$ is $K_0-$bilipschitz, for some constant $K_0$ depending on
the pair $(K,C)$.

Let $H$ be the
cyclic group generated by the matrix $\left(
\begin{array}{cc} p & 0 \\ 0 & \frac{1}{p} \end{array} \right)$.  
We make the following definitions relating to a group-invariant
diameter function which will allow us to state Theorem
\ref{thm:action:rigidity} (Action Rigidity).

\medskip

\noindent
{\bf Definitions.}

\bigskip
\noindent
1. For any subset $S$ of $ \R \times \Qp$ and $H$ as above,
define the $H$-invariant diameter of $S$ by $$\delta_H(S) = \inf_{T \in H} diam
(T(S)).$$

\noindent
For the remainder of this paper we will write $\delta$ for $\delta_H$.

\bigskip
\noindent
2. For subsets $S_1,S_2$ of $ \R \times \Qp$, we say that the map $\phi:
S_1 \rightarrow S_2$ is {\it quasi-adapted} to $\delta$ if there exists a
map $\alpha: {\bf N} \rightarrow {\bf N}$
such that for any compact set $V$, 
$$\delta(V) \leq k \Rightarrow \delta(\phi(V)) \leq \alpha(k)$$
and
$$\delta(\phi(V)) \leq k \Rightarrow \delta(V) \leq \alpha(k).$$

\bigskip
\noindent
3. Let $S$ be a subset of $\dq$.  We say that $S$ has
  {\it bounded height} if $S \subset \cz$ for some $M \in \Z^+$ with $(M,p) = 1$.  

\bigskip
\noindent
4. A bijection $\phi: \dq \rightarrow \dq$ is said to be
{\it quasi-integral} if both $\phi$ and $\phi^{-1}$ take sets of bounded
height to sets of bounded height.  

\bigskip
\noindent
5. A bijection $\phi: \dq \rightarrow \dq$ is said to be {\it quasicompatible} with
$H$ if $\phi$ is quasi-integral and, when restricted to sets of bounded
height, both $\phi$ and $\phi^{-1}$ are quasi-adapted to $\delta$.

\bigskip

In a sequence of lemmas we show that the bilipschitz map $\phi: \dq \rightarrow \dq$ obtained from the original
quasi-isometry $f$ is quasi-compatible with $H$.

Associate to each point of $\dq$ the horosphere boundary component
of $\Omega_p$ based at that point.  
Then the action of $\Gamma_{\infty}$ on $\partial \Omega_p$ induces an action
on $\dq$.  
In particular, this action preserves denominators,
i.e. $\Gamma_{\infty} \cdot \frac{1}{M} \dz \subset \frac{1}{M} \dz$.
Consider the diameter function on subsets $S$ of $\dq$ defined by 
$$\delta_{\Gamma}(S) = \inf_{T \in \Gamma_{\infty}} diam \ 
  T(S).$$  
Also consider the following diameter function, with $S$ as above.  
Let $\delta_B(S) $ be the diameter in $\Omega_p$ of the smallest
metric ball which intersects all horosphere boundary components based
at points of $S$.  
We now show that these two diameter functions are {\em quasi-identical} when
restricted to $\frac{1}{M} \dz$, i.e. if $\delta_{\Gamma}(S)$ is small
for a  subset $S$ of $\dq$, then $\delta_B(S)$ is bounded and vice
versa.  

\begin{lemma}
The restrictions of $\delta_{\Gamma}$ and $\delta_B$ to $\frac{1}{M}
\dz$ are quasi-identical.
\end{lemma}

\begin{proof}
Let $S$ be a subset of $\dq$ so that $\delta_B(S)$ is small, say
$\delta_B(S) = \epsilon$.  
Then there exists a point $x \in \Omega_p$ which is within $\epsilon$
of all horosphere boundary components based at points of $S$.  
Since $\Omega_p / \PZp$ is compact by construction, there are only a
finite number of choices for the set $S$ (modulo Aut($\dz$)).  
Thus $\delta_{\Gamma}(S)$ must be bounded.

Now suppose $\delta_{\Gamma}(S)$ is small.  By compactness there are
only finitely many choices for $S$ (modulo Aut($\dz$)).  Hence
$\delta_B(S)$ is bounded.  
\end{proof}


\begin{lemma}
\label{lemma:quasi:integral}
For every $M \in \Z^+$, there exists  $M' \in Z^+$ depending on
the bilipschitz constant $K_0$ of $\phi$ and the space $\Omega_p$, so that
$\phi(\frac{1}{M} \dz) \subset \frac{1}{M'} \dz$.  
\end{lemma}

\begin{proof}
Consider the action of $\Gamma_{\infty} = Stab_{\PZp}(\infty)$ on
$\frac{1}{M} \Delta$.  
This action preserves denominators, i.e. $\Gamma_{\infty} \cdot
\frac{1}{M} \Delta \subset \frac{1}{M} \Delta$.  
In particular, under this action $\frac{1}{M} \Delta
/_{\Gamma_{\infty}}$ consists of a finite set of points, which we view
as a finite collection of horospheres $\sigma_1, \cdots \sigma_n
\subset \partial \Omega_p$.  
Choose a point $x \in \si$ and consider the smallest metric ball in
$\Omega_p$ based
at $x$ which intersects all of the $\sigma_i$.  
Let $\epsilon $ be the radius of this ball.  
There is an $\epsilon'$ (depending on $\epsilon$) so that the
$\epsilon'$ ball around $f(x)$ must intersect all of the
$f(\sigma_i)$.  
Thus there are only a finite number of choices of horospheres in
$\partial \Omega_p$ for the images $f(\sigma_i)$. 
It follows that there is a number $M' \in \N$ so that the collection of horospheres
$\{ f(\sigma_i) \}$ is based at points in $\frac{1}{M'} \dz$.  
Since $\Gamma_{\infty}$ preserves denominators, the image of
$\sigma_{\alpha}$ (for any $\alpha \in \frac{1}{c} \dz$) must be a
horosphere based at a point of $\frac{1}{M'} \dz$.  
This is equivalent to saying that $\phi(\frac{1}{M} \dz) \subset
\frac{1}{M'} \dz$.  
\end{proof}

Lemma \ref{lemma:quasi:integral} shows that $\phi$ is quasi-integral, and since $\phi$ is a
bilipschitz map of both $\R$ and $\Qp$, it is quasi-adapted to
$\delta$.  
Hence $\phi$ is quasicompatible with $H$.

We can now state Theorem \ref{thm:action:rigidity}.
We say that a map $\psi: \R \times \Qp \rightarrow \R \times \Qp$ is
affine if its restriction to each factor is affine.  This theorem is
an $2-$dimensional $S$-arithmetic version of the Action Rigidity Theorem of
R. Schwartz \cite{S}.

\begin{theorem}[Action Rigidity]
\label{thm:action:rigidity}
Let $\dq \subset \R \times \Qp$ be the lattice $\{ (a,a) | a \in \Q
\}$ and let $H$ be the group
generated by the matrix  $  \left(
\begin{array}{cc} p & 0 \\ 0 & \frac{1}{p} \end{array} \right)$.
Then any bilipschitz map $\phi: \dq \rightarrow \dq$ which is
quasi-compatible with $H$ is the restriction of an affine
  map of $\R \times \Qp$.
\end{theorem}

\subsection{The Parallelogram Lemma}
\label{sec:plemma}

The key lemma in the proof of Theorem \ref{thm:action:rigidity} is
Lemma \ref{lemma:parallelogram} (Parallelogram Lemma). 
Let $\dq$ be as in \S \ref{sec:background:ar} and let $\phi: \dq
\rightarrow \dq$ be quasi-compatible with $H$ and $K_0-$bilipschitz.  
We now make some preliminary definitions.  

\noindent
{\bf Definition.} Let $M \in \Z^+$.  A {\em parallelogram} $P$ in $\cz$ is a quadruple of points $\left[
\begin{array}{cc} a & b \\ c & d \end{array} \right]$ with $a,b,c,d
  \in \cz$ satisfying $a -c = b -d$.
\bigskip

We say that $\phi(P)$ is the quadruple $\left[
\begin{array}{cc} \phi(a) & \phi(b) \\ \phi(c) & \phi(d) \end{array}
\right]$.  
The goal of Lemma \ref{lemma:parallelogram} is to determine when $\phi(P)$ is
itself a parallelogram.  
We now define two quantities associated to a
parallelogram $P$ which will be invariant under the group action and
translation.  

\bigskip
\noindent
{\bf Definition.} Let $P = \left[
\begin{array}{cc} a & b \\ c & d \end{array} \right]$ be a
  parallelogram.  
The {\em $H$-invariant perimeter} of $P$ is given by
$$ per(P) = \delta(a \cup b) + \delta(a \cup c).$$  
The {\em shape} of $P$ is given by $$ shape(P) =  |\nu(b-a) 
- \nu(c-a)|$$
where $\nu(p^n \frac{x}{y}) = n$ is the $p-$adic valuation on $\Qp$.  
For any $T \in H$ and $x \in \dq$ we have $per(P) = per(T(P) + x)$ and
$shape(P) = shape(T(P) + x).$

\bigskip
\bigskip

Let $P = \left[
\begin{array}{cc} a & b \\ c & d \end{array} \right]$ be a
  parallelogram such that $\phi(P) = \left[
\begin{array}{cc} a' & b' \\ c' & d' \end{array} \right]$ is again a parallelogram.  
Since $\phi$ is quasi-compatible with $H$, there exists a constant $D$ (depending on
$per(P)$ and $K$) such that 
$$ \delta(a' \cup b') + \delta(a' \cup c') \leq D.$$
By symmetry we also have $\delta (c' \cup d') + \delta ( b' \cup d') \leq D$.  
From lemma \ref{lemma:quasi:integral} we obtain a constant $k$ such
that $\phi(\cz) \subset \kz$.

We now describe the points $x \in \kz, \ (k,p) = 1$, satisfying
$\delta(0 \cup x) \leq D$.  
$$S_{k,D} = \{ x \in \kz | x \neq 0 \ and \ \delta(0 \cup x ) \leq D\}.$$
The set $S_{k,D}$ is the orbit under $G$ of a finite set of points of
$\kz$, denoted $\{\frac{a_i}{kp^{r_i}} \}_{i \in I}$, where $(a_i, p) = 1$.
Since $\delta$ is invariant under translation, the points of $\kz$
within $D$ of some point $y \in \kz$ are given by $S_{k,D} + y$.  


We use the following notation in the statement of Lemma
\ref{lemma:parallelogram}.  
Let 
$$B_1 = max_{i,j \in I} \{a_i - a_j\}, \ B_2
= max_{i \in I} \{a_i\},
$$
$$ \ B_3 = max_{i,j \in I}  \{ |\nu (a_i - a_j )|
, |\nu( a_i + a_j)| \}$$
and 
$$\ B_4 = max_{i,j,k \in I} \{ \nu(a_i+a_j-a_k) \}.$$

\begin{lemma}[Parallelogram Lemma]
\label{lemma:parallelogram}
Let $\phi: \dq \rightarrow \dq$ be quasi-compatible with $H$ and $K_0-$bilipschitz.  
 Let $\log _p (K_0) = R$ and 
$L \in {\bf N}, \ (L,p) = 1$.  If $P$ is a parallelogram in $\frac{1}{L}
\dz$ with $$per(P) \leq L \  and \ shape(P) > s_0$$
where $$s_0  = max \{ 2\log_p(B_1 + B_2) + 2R, \
3B_3^2 + 2R, \ 2B_4 + 2R \}$$
then $\phi (P)$ is also a parallelogram.
\end{lemma}

\begin{proof}
Let $P = \left[
\begin{array}{cc} a & b \\ c & d \end{array} \right] $ and $\phi(P) = \left[
\begin{array}{cc} a' & b' \\ c' & d' \end{array} \right] $.  
We are assuming that $\phi(0) = 0$, so without loss of generality
translate $P$ so that $a=0$ and hence $a' = 0$.  
As stated above, there is a constant $D$ so that
$\delta(0 \cup b') + \delta(0 \cup c') \leq D$ and $\delta (b' \cup
d') + \delta (c' \cup d') \leq D$.

We can write $b' = p^n \frac{x}{k}$ and $c' = p^m \frac{y}{k}$ where $x=a_i$ and
$y=a_j$ for some $a_i,a_j$ as above.  
Then $d'$ can be expressed in one of two ways.  
It is $\delta-$close to both $b'$ and $c'$, hence of the form $p^n
\frac{x}{k} + p^N \frac{z}{k}$ and also $p^m \frac{y}{k} + p^M \frac{w}{k}$ where $z=a_i, \ w=a_j$ for some $i,j \in I$.  
Setting these expressions equal and clearing denominators yields
$$p^nx + p^N z = p^my + p^M w. \ \ \ \ \ \ \ \ \
\ (*)$$
Note that $x,y,z,w$ are all relatively prime
to 
$p$.  

We first obtain a lower bound on $n-m$, with $n,m$ as above.  
We can write $b = p^{n_0} \frac{x_0}{L}$ and $c = p^{m_0}
\frac{y_0}{L}$, with $(x_0,p) = 1$ and $(y_0, p) = 1$, chosen so that
$shape(P) = n_0-m_0$.  
We know by assumption that $n_0 - m_0 > s_0$.  
Since $\phi$ is a $K_0-$bilipshitz map on $\Q \subset \Qp$ and
$\phi(0) = 0$, letting $R = \log _p(K_0)$ we obtain
$$p^{-n_0-R} \leq | \phi( p^{n_0} \frac{x_0}{L})| _p = |p^n
\frac{x}{k}|_p \leq p^{-n_0 + R}.$$
Hence $n_0-R \leq n \leq n_0 + R$; similarly, $m_0-R \leq m \leq m_0
+ R$.  
Thus, $$ (n_0 - m_0) - 2R \leq
n-m \leq (n_0 - m_0) + 2R.$$
Using the fact that $n_0-m_0 > s_0$, we see that
$$n-m \geq s = max \{2\log_p(B_1 + B_2),\ 3B_3^2, \ 2B_4 \}.$$

Consider again the expression $(*)$ for $d'$.  First we
assume that not all exponents are equal.  
If all the exponents were equal, then $n$ would
equal $m$, which is impossible since $n-m>s>0$.  
Suppose first that $n>N$.  
Since the $p-$norms of both sides must be equal, we know that $N=min
\{ m, M \}$.  

\medskip
\noindent
{\bf Case 1.} Suppose $N=m$ and $M \neq m$.  Then $(*)$ simplifies to $$p^{n-m}x =
 p^{M-m}w + (y-z).$$  
We know that $n-m > s$. 
Suppose $y-z \neq 0$.  
The highest power of $p$ dividing the right hand side of the
equation is bounded by $B_3$.  
But by the choice of $s$, we have $n-m > B_3$; hence this situation cannot occur.

If by chance $M-m = B_3$, and $y-z = p^{B_3}t$ for some integer $t$,
then note that $\nu(w+t)$ is also bounded.  If this is the case,
change the initial value of $s_0$ by the quantity $\nu(w+t)$, and we may assume
that we are not in this case.  

So we must have $y-z = 0$, i.e. $z=y$ and $p^{n-m}x = p^{M-m}w$.  
Thus $n=M$ and $x=w$, and $\phi(P)$ is a parallelogram.

\medskip
\noindent
{\bf Case 2.}  Suppose $N=M$ and $M \neq m$.  
Then $(*)$ simplifies to $$p^{m-M}(p^{n-m}x-y) = w-z.$$  
Since $m-M > 0$
and $w-z \leq B_1$, 
our choice of $s_0$ insures that $p^{m-M}(p^{n-m}x-y) > B_1$, a
contradiction.  

\medskip
\noindent
{\bf Case 3.}  Suppose that $N=M=m$. 
Then $(*)$ simplifies to $$p^{n-m}x = y+z-w.$$  
By our choice of $s_0$, the exponent $m-n > B_4$, a contradiction.

Now suppose that $n<N$.  
Then we have $n = min \{m,M \}$.  
Since $n-m > s>0$, we cannot have $n=m$.

\medskip
\noindent
{\bf Case 4.}  Suppose that $n=M$.  
But then $M < m < n=M $, a contradiction.

Lastly, we consider the case $n = N$.  

\medskip
\noindent
{\bf Case 5.}  Suppose $n = N$.  
Then $(*)$ simplifies to $$p^n(x+z) = p^my + p^M w.$$  
The quantity $x+z$ can take one of a finite number of values, so we
can write $x+z = p^h f$ where $(f,p) = 1$ and $h$ is bounded by
$B_3$.  
Then we must have $n+h = min(m,M)$.  
If $n+h = m$, then $-h = n-m > s$ which cannot happen, since $|h| \leq
B_3$ but $n-m > B_3$.  
If $n+h = M$, then $(*)$ becomes $p^{n-m+h} (f-w) = y$.  
Again by the choice of $s_0$, the left hand
side of the equation is greater than the right hand side.
\end{proof}

\subsection{Proof of the Action Rigidity Theorem}
\label{sec:actionrigidity}

We now prove Theorem \ref{thm:action:rigidity}.

\bigskip
\noindent
{\it Proof of Theorem \ref{thm:action:rigidity}.} Fix $q \in \N, \ (q,p) = 1$ and let $S$ be a generating set for $\qz$
containing $\frac{1}{q}$. 
Let $s_0$ be as in the Lemma \ref{lemma:parallelogram}, and let $H(S)$ denote the orbit of $S$ under $H$.  
Given $x,y \in \qz$ we say that $(x,y)$ is a distinguished pair if
$x-y \in H(S)$.  
Let $C = max_{a \in H(S)} (q,2\delta(0 \cup a))$.  

We write $P((x,y),(z,w))$ if $P = \left[ \begin{array}{cc} x & y \\ z
& w \end{array} \right]$ is a parallelogram with $P \subset \qz$,
$per(P) \leq C$ and $shape(P) \geq s_0$.  
Then by the Parallelogram lemma, $\phi(P)$ will also be a
parallelogram, i.e. $$\phi(x) - \phi(z) =
\phi(y) - \phi(w).$$

We write $\overline P((x,y),(z,w))$ if there is a finite sequence of
pairs $(a_i,b_i)$ for $i = 0, \cdots, n$ with $(a_0,b_0)
= (x,y)$ and $(a_n,b_n) = (z,w)$ such that
$P((a_i,b_i),(a_{i+1},b_{i+1}))$.  
This means that we have a sequence of parallelograms between the two
given pairs of points, each satisfying the conclusions of the
parallelogram lemma, and with any two
consecutive parallelograms in this sequence sharing a common side.
Concatenating these intermediate parallelograms allows us to conclude
that the parallelogram given by the original pairs
also satisfies the conclusions of the Lemma \ref{lemma:parallelogram}.

We first need the following lemma.

\begin{lemma}
Let $a \in \qz$ be arbitrary, and let $(u,v)$ be a distinguished pair.
Then $\overline P((u,v),(u+a,v+a))$.  
\end{lemma}

\begin{proof}
Consider $P_y = \left[ \begin{array}{cc} u & v \\ u + y & v + y
\end{array} \right]$ for any $y \in H(S)$.  
By construction, we have $per(P_y) \leq C.$
Let $Y=Y(u,v) \subset H(S)$ denote those $y \in H(S)$ such that
$shape(P_y) \geq s_0$.  
Let $\Sigma_0Y$ be the sublattice in $\dq$ generated by $Y$ and
$\Sigma Y = \Sigma_0Y \cap \qz$.  
We will first show that $\overline P ((u,v),(u+x,v+x))$ for any $x \in
\Sigma Y$ and then show that $\Sigma Y = \qz$.  
We prove this first assertion by induction.  If $x \in Y$, it is certainly true that $P((u,v),(u+x,v+x))$.  
If $x \in \Sigma Y$, write $x=x'+y$ for $y \in Y$, 
and we have $\overline
P((u,v),(u+x',v+x'))$ by the induction hypothesis.  Consider the parallelogram 
$$P = \left[ \begin{array}{cc} u+x' & v+x' \\ u + x' +y & v + x' +y
\end{array} \right] = \left[ \begin{array}{cc} u+x' & v+x' \\ u + x & v + x
\end{array} \right].$$
We see that $shape(P) = shape(P_y) \geq s_0$.  
Also, $\delta (u+x' \cup u + x) = \delta (0 \cup y) \leq \frac{1}{2}
C$ and $\delta (u+x' \cup v + x') = \delta (u \cup v) \leq \frac{1}{2}
C$.  
So we have $per(P) \leq C$.  
Therefore we know that $P((u+x',v+x'),(u+x,v+x))$.  

Now we show that $\Sigma Y = \qz$.  
We know that $shape(P_y) = |\nu(u-v) - \nu(y)|$.  Let $\omega =
\nu(u-v)$.  
Then $Y= \{ y \in H(S) | |\omega - \nu(y)| \geq s_0 \}$.  
Hence for large enough $n$, we have $\frac{1}{p^n} \frac{1}{q} \in Y$, so $\Sigma Y
= \qz$, proving the lemma.
\end{proof}

Given $x \in \qz$, we can create a distinguished pair by taking
$(x,x+a)$, for any $a \in S$.  
So for any $x,y \in \qz$ we know that $P((x,x + a),(y,y +
a))$, or $\phi(x + a) - \phi(x) = \phi(y +
a)- \phi(y)$.  
Since $S$ generates $\qz$ we know that for any $x,y,z \in
\qz$ we have $\overline P((x,x+z),(y,y+z))$ or 
$$\phi(x+z) - \phi(x) =
\phi(y+z) - \phi(y).$$  
In particular, since $\phi(0) = 0$, we know that $\phi|_{\qz}$ is
multiplication by a constant $C_q$.  
Let $q_1, q_2 \in \N$, with $(q_1, p) = (q_2, p) = 1$.  
Then we must have $C_{q_1q_2} = C_{q_1} = C_{q_2}$.  Hence there is a
constant $\alpha$ such that the map $\phi|_{\dq}$ is multiplication by
$\alpha$.  
It follows that $\phi$ is the restriction of an affine map of $\R
\times \Qp$.  
\qed

\subsection{Proof of Theorem A}
\label{sec:proof:mainthm}

We now use Theorem \ref{thm:action:rigidity} to prove Theorem A.

\bigskip
\noindent
{\it Proof of Theorem A.} It is clear that every commensurator of $\PZp$ gives rise to a unique
quasi-isometry of $\Omega_p$.  
Given $f \in QI(\PZp)$, we will now choose a commensurator 
$g \in Comm(\PZp)$ so that the composition $g \circ f$
is a bounded 
distance from the identity map.  

Let $f \in QI(\PZp)$.  
Then $f$ induces a quasi-isometry from $\Omega_p$ to itself.    
Compose $f$ with a commensurator so that $f(\si) = \si$ and $f(\so) = \so$.   Combining Theorems
\ref{thm:BS:2} and \ref{thm:action:rigidity}, we know that $f$
induces a map $f_l$ on the lower boundary of $\si$ which is
multiplication by a constant $\alpha$ and determines the permutation
of the horospheres under $f$. 

Recall from \S \ref{sec:algebraic} that we view $PSL_2(\Q)$ as
the $\Q-$points $G'_{\Q}$, where $G' = Ad(SL_2(\Cp))$.  
Compose $f$ with the commensurator $g = Ad\left( \begin{array}{cc}
\frac{1}{\sqrt{\alpha}} & 0 \\ 0 & \sqrt{\alpha} \end{array} \right)$.  
Then the permutation of the horospheres under $g \circ f$ is the
identity permutation.  
We must now show that for any $m \in \Omega_p$, the image $f(m)$ lies in
an $\epsilon$-ball around $m$, where $\epsilon$ is independent of the
choice of $m$.  

The set of points within $n$ units of any three distinct horospheres
has bounded diameter, independent of the choice of horospheres.  
Also, any quasi-isometry $f$ has the property that 
$$
hd(f(Nbhd(A) \cap Nhbd(B)), Nbhd(f(A) \cap Nbhd(B))) < \infty
$$  
where $hd$ denotes Hausdorff distance.  
So if $x \in \Omega_p$ lies in the intersection of the
$n$-neighborhoods of three horospheres, then there is a constant $n'$
so that $f(x) \in B_{n'}(x)$.  Thus $f \circ g$ is a bounded distance
from the identity map, so the natural map $\Psi : Comm(\PZp)
\rightarrow QI(\PZp)$ is an isomorphism.  
\qed

\subsection{Proof of Theorem C}
\label{sec:pf:C}

We now prove Theorem C.  
The proof uses some standard techniques from the study of
quasi-isometric rigidity for lattices in semisimple Lie
groups.  In addition to these techniques, Theorem A is applied, as well
as the $S$-arithmetic superrigidity theorem of Margulis \cite{M}.    

\bigskip
\noindent
{\it Proof of Theorem C.}  
Since $\PZp$ is quasi-isometric to the space $\Omega_p$, we get
a quasi-isometry $f: \Gamma \rightarrow \Omega_p$.  
To obtain the exact sequence of the theorem, we will find a representation
$\rho: \Gamma \rightarrow \G$ with finite kernel so that
$\rho(\Gamma)$ is a nonuniform lattice in $\G$, hence commensurable
to $\PZp$ by \cite{M}.  

Since $\Gamma$ acts on itself by isometries via left multiplication
$L_{\gamma}$, for all $\gamma \in \Gamma$, we obtain a uniform family
of quasi-isometries
$$f_{\gamma} = f \circ L_{\gamma} \circ f^{-1} : \Omega_p \rightarrow
\Omega_p.$$
From Theorem A, we can think of $f_{\gamma}$ as a commensurator,
hence a bounded distance from an isometry.  
For each $f_{\gamma}$, compose with a bounded alteration $B_{\gamma}$
to obtain a map
$$\rho: \Gamma \rightarrow Isom(\HT).$$
First we show that $\rho$ is a homomorphism.  
Suppose $\rho(\gamma)$ is a bounded distance from the identity
isometry, i.e. $d_{\HT}(x, \rho(\gamma) \cdot x) < \epsilon$ for all $x
\in \HT$.  
But then it follows that $\rho(\gamma)$ is the identity, hence $\rho$ is
a homomorphism.  
We must show that $\rho(\gamma)$ is a lattice and that $\rho$ has
finite kernel.

Since $f(\Gamma)$ is a net in $\Omega_p$, $\rho(\Gamma) = \{ f_{\gamma}
\} = \{ B_{\gamma} \circ f \circ L_{\gamma} \circ f^{-1} | \gamma \in \Gamma \}$ acts
cocompactly on $\Omega_p$.  
It follows that $\rho(\Gamma)$ acts on $\HT$ with cofinite volume.  

Choose a basepoint $x \in \Omega_p$ and consider $d_{\HT}(x,
\rho(\gamma) \cdot x)$.  
For finitely many $\gamma \in \Gamma$, $\rho(\gamma)$ moves $x$ a
bounded amount, i.e. there exists a constant $C' > 0$ so that
$d_{\HT}(x, \rho(\gamma) \cdot x) \leq C'$ for some finite set $\{
\gamma_i \} \subset \Gamma$.  
(To find $C'$ we are using the fact that $\{ f_{\gamma} \}$ is a
uniform family of quasi-isometries.)  
In particular, $d_{\HT}(x, \rho(\gamma) \cdot x) = 0$ for only
finitely many $\gamma \in \Gamma$.  
But if $d_{\HT}(x, \rho(\gamma) \cdot x) > 0$, then $\rho(\gamma)$ is
not the identity isometry.  
Hence $\rho$ has finite kernel.  The previous paragraph showed that
$\rho(\Gamma)$ is discrete, so we can now construct the short exact
sequence
$$1 \rightarrow N \rightarrow \Gamma \rightarrow \Lambda \rightarrow
1$$
where $N = ker(\rho)$ and $\Lambda = \rho(\Gamma)$ is a lattice in
$\G$.

To complete the proof, we apply the following theorem of Margulis
\cite{M} (Theorem \ref{thm:superrigidity}), stated in the case $n = 2$, from which we conclude that the lattice $\Lambda$
obtained above is commensurable to $\PZp$.  

\begin{theorem}\cite{M}
\label{thm:superrigidity}
Let $p$ be a prime and $n \in N^+, \ n \geq 2$.  If a subgroup
$\Gamma$ of $SL_2(\Q)$ under the diagonal embedding in $SL_2(\R)
\times SL_2(\Qp)$ is a lattice in $SL_2(\R)
\times SL_2(\Qp)$, then the subgroups $\Gamma$ and $SL_2(\F)$ are
commensurable.  
\end{theorem}

\bigskip
\noindent
The proof of Theorem \ref{thm:superrigidity} uses Margulis'
superrigidity theorem for $S$-arithmetic groups.
Theorem \ref{thm:superrigidity} completes the proof of Theorem C.
\qed

\section*{APPENDIX}

We now compute the commensurator group $Comm(\PZp)$ of $\PZp$ in $\G$.

\noindent
{\bf Proposition.}
\label{prop:comm:group}
{\it The commensurator subgroup for $\PZp$ is given by 
$$Comm(\PZp) = PSL_2(\Q) $$ 
$$= \{ (a,a) | a \in PSL_2(\Q) \} \subset
PSL_2(\R) \times PSL_2(\Qp).$$}

\begin{proof}
Since any commensurator $g$ must preserve $\Omega_p$, hence preserve the
horospheres, i.e. map $\Q \cup \{ \infty \}$ to $\Q \cup \{ \infty \}$, we must have $g \in PSL_2(\Q)$.  
We now show that the commensurator subgroup, which is contained in
$PSL_2(\R) \times PSL_2(\Q_p)$, is exactly the set $\{ (g,g) | g \in
PSL_2(\Q) \}$.  

First we show that for $b \in PSL_2(\Q)$, the element $h = (Id, b)$ is not a
commensurator of $\PZp$.  
We know that an element $(a,b) \in Comm(\PZp)$ must induce a
quasi-isometry on the space $\Omega_p$.  
Since $a \in PSL_2(\Q) \subset PSL_2(\R)$ and $b \in PSL_2(\Q) \subset
PSL_2(\Qp)$, we
can let $a$ act by conjugation on the $\A$ factor and $b$ by
conjugation on $T_p$.  
So $h$ acts by the identity on $\A$, hence $h(\sa) = \sa$, although
$\sa$ is not fixed pointwise.  
In addition, $h(\si) = \si$, so the closeness line between $\sa$ and
$\si$ is also preserved.  
Since this is true for all $\alpha \in \Q$, if we view $h$ as a map on
$\Qp$, then $h$ fixes a copy of $\Q \subset \Qp$.  
Hence when we extend $h$ continuously to $\Qp$, we see that it must be
the identity.

We now show that if $a \in PSL_2(\Q)$, then $(a,a)$ is a commensurator of
$\PZp$.  
A calculation shows
that for any element 
$$(g,g) \in \PZp \times \PZp \subset \G$$ the
expression $$(a,a)(g,g) (a^{-1},a^{-1})$$ can be written
as a matrix all of whose entries have bounded denominator, say bounded
by $d$, meaning that each denominator is of the form $dp^r$ for some
$r \in \Z^+$.  
Thus the subgroup $$(a,a)(\PZp, \PZp) (a^{-1},a^{-1})$$ is of finite
index in $\PZp \times \PZp$.  
Hence $(a,a)$ is a  commensurator.  

Now we will show that $ Comm( \PZp )$ consists entirely of elements of
the form $(a,a)$, with $a \in PSL_2(\Q)$.  
Suppose $(a,b) \in Comm(\PZp)$, with $a,b \in PSL_2(\Q)$.  
We know that $(a^{-1},a^{-1})$ is an element of $ Comm(\PZp)$.  Compose these two
elements to form a new commensurator $(a,b) \circ (a^{-1},a^{-1}) =
(Id, ba^{-1})$.  
By the reasoning above, if this is a commensurator, then we must have
$ba^{-1} = Id$, which implies that $b=a$, so our original commensurator
must have been of the form $(a,a)$ with $a \in PSL_2(\Q)$.
\end{proof}

\noindent
Jennifer Taback\\
Dept. of Mathematics\\
University of California-Berkeley\\
Berkeley, CA 94720\\
E-mail: jen@math.berkeley.edu


\begin{thebibliography}{FM2}





\bibitem[F]{F}
B. Farb, The quasi-isometry classification of lattices in semisimple
Lie groups, to appear in {\it Math. Res. Letters}.

\bibitem[FM]{FM}
B. Farb and L. Mosher (appendix by D. Cooper), A rigidity theorem for
the solvable Baumslag-Solitar groups, {\it Inventiones Math},
Vol. 131, No. 2 (1998), pp. 419-451.

\bibitem[FM2]{FM2}
B. Farb and L. Mosher, Quasi-isometric rigidity for the solvable
Baumslag-Solitar groups, II, preprint. 

\bibitem[FS]{FS}
B. Farb and R. Schwartz, The large-scale geometry of Hilbert modular
groups, {\it J. Diff. Geom.} 44, No. 3 (1996), pp. 435-478.

\bibitem[G]{G}
M. Gromov, Infinite groups as geometric objects, Plenary Address, {\it
  Proc. of the I.C.M.}, Warsaw, 1983.

\bibitem[LMR]{LMR}
A. Lubotzky, S. Mozes, M.S. Raghunathan, Cyclic subgroups of
exponential growth and metrics on discrete groups, {\it
  C.R. Acad. Sci. Paris - Ser. I Math.} 317 no 8, 1993.

\bibitem[M]{M}
G.A. Margulis, {\it Discrete Subgroups of Semisimple Lie Groups,}
Springer-Verlag, 1991.

\bibitem[Mi]{Mi}
J. Milnor, A note on curvature and fundamental group,  {\it
  J. Diff. Geom.} Vol. 2 (1968), pp. 1-7.

\bibitem[S1]{S}
R. Schwartz, Quasi-isometric rigidity and Diophantine approximation,
{\it Acta Mathematica} Vol. 177 (1996), pp. 75-112.

\bibitem[S2]{S2}
R. Schwartz, The quasi-isometry classification of rank $1$ lattices,
{\it IHES Sci. Publ. Math.,} Vol. 82, 1996.


\bibitem[Se]{Se}
J.P. Serre, {\em Trees}, translated by J. Stillwell, Springer-Verlag, 1980.



\end{thebibliography}
\end{document}